\documentclass[9points, conference]{IEEEtran}
\usepackage[latin9]{inputenc}
\usepackage[active]{srcltx}
\usepackage{color}
\usepackage{float}
\usepackage{bm}
\usepackage{amsmath}
\usepackage{amssymb}
\usepackage{graphicx}
\usepackage{amsfonts}
\usepackage{hyperref}
\usepackage{epstopdf} 
\usepackage{multirow}

\usepackage{tabularx}
\usepackage{tcolorbox}

\let\labelindent\relax
\usepackage{enumitem}

\makeatletter

\newcommand{\norm}[1]{\left\lVert#1\right\rVert}

\usepackage{xcolor}
\usepackage{colortbl}
\usepackage{lipsum}

\floatstyle{ruled}
\newfloat{algorithm}{tbp}{loa}
\providecommand{\algorithmname}{Algorithm}
\floatname{algorithm}{\protect\algorithmname}

\makeatletter
\let\oldfootnote\footnote
\def\footnote{\@ifstar\footnote@star\footnote@nostar}
\def\footnote@star#1{{\let\thefootnote\relax\footnotetext{#1}}}
\def\footnote@nostar{\oldfootnote}
\makeatother

\usepackage[noadjust]{cite}

\usepackage[T1]{fontenc}
\usepackage{float}
\usepackage{subfigure}
\makeatletter

\floatstyle{ruled}
\newfloat{algorithm}{tbp}{loa}
\providecommand{\algorithmname}{Algorithm}
\floatname{algorithm}{\protect\algorithmname}

\newtheorem{theorem}{Theorem}

\newcounter{algo}

\ifCLASSINFOpdf
\else
\fi
\usepackage[footnotesize]{caption}
 \usepackage{flushend}

\makeatother

\begin{document}

\title{\vspace{-0.2cm} Distributed Dictionary Learning\vspace{-0.2cm}}

\author{Amir Daneshmand, Gesualdo Scutari, and  Francisco Facchinei$^\dagger$\vspace{-0.3cm}%

} 
 
\maketitle
\footnote*{$^\dagger$Daneshmand and Scutari are with the School of Industrial Engineering, Purdue University, West-Lafayette, IN, USA; emails: \texttt{<adaneshm, gscutari>@purdue.edu}. Facchinei is with the Dept. of Computer, Control, and Management Engineering, University of Rome ``La Sapienza'',  Rome, Italy; email: \texttt{francisco.facchinei@uniroma1.it}.   The work of  Scutari was supported by the USA National Science Foundation under Grants CIF 1564044 and CAREER Award 1555850, and the ONR N00014-16-1-2244. 
}
\begin{abstract}
The paper studies \emph{distributed Dictionary Learning (DL)} problems where 
the    learning task is distributed over a multi-agent network with time-varying (nonsymmetric) connectivity. This formulation is relevant, for instance,  in
Big Data scenarios where  massive amounts of data are   collected/stored
in different spatial locations and it is unfeasible to aggregate and/or process all data in a fusion center, due to resource limitations, communication overhead or privacy considerations.
We develop a general \emph{distributed} algorithmic framework for the (nonconvex) DL problem and establish its asymptotic convergence.   The new  method hinges on Successive Convex Approximation (SCA) techniques coupled with i) a gradient tracking mechanism  instrumental
to locally estimate the missing global information; and ii) a consensus step, as a mechanism to distribute the computations   among the agents. To the best of our knowledge, this is the first distributed
algorithm with provable convergence for the DL problem and, more in
general, bi-convex optimization problems over (time-varying) directed graphs. 
 \end{abstract}

\IEEEpeerreviewmaketitle
\section{Introduction\label{sec:intro} }
\IEEEPARstart{T}{he} dictionary learning problem \cite{Elad_book10} consists in finding  
an overcomplete basis (a.k.a the dictionary)   by which a given set of data 
can be sparsely represented. 
This technique can be leveraged to solve many machine learning and inference tasks, including  image denoising/debluring/inpainting, superresolution \cite{Elad_Aharon_2006,Mairal_Bach_2010},  dimensionality reduction
\cite{Shen_Huang08}, bi-clustering \cite{Lee_Shen_Huang_Marron10}, feature-extraction and classification
\cite{Mairal_ICASSP08}, and prediction \cite{NIPS2008_3538}.
 
  Denoting by $\mathbf{S}\in \mathbb{R}^{M\times N}$ the data (observation) matrix, by $\mathbf{D}\in \mathbb{R}^{M\times K}$ the overcomplete dictionary matrix, and by $\mathbf{X}\in \mathbb{R}^{K\times N}$ the sparse representation matrix of the signal, the DL problem for sparse representation reads  (see, e.g., \cite{elasticNet_Hastie_2005})\vspace{-3pt}
\begin{equation}
\hspace{-0.3cm}\begin{array}{rl}\displaystyle
\min_{\mathbf{D},\mathbf{X}}\quad& \dfrac{1}{2}\norm{\mathbf{S}-\mathbf{D}\mathbf{X}}^2_F+\lambda\norm{\mathbf{X}}_{1,1}+ {\mu} \norm{\mathbf{X}}_F^2
\\[0.5em]
\mathrm{s.t.}\quad &\mathbf{D} \in {\cal D} \triangleq
\{\mathbf{D}\,:\,\norm{\mathbf{D}\mathbf{e}_k}_2\leq\alpha, k=1,2,\ldots,K\},
\tag{P1}\label{eq:P1}
\end{array} 
\end{equation}
where  $\norm{\mathbf{X}}_{1,1}\triangleq \sum_{i,j}|X_{i,j}|$, and
$\mathbf{D}$ is constrained to belong to the set $\mathcal D$  to avoid unbound solutions, with  $\alpha>0$ and $\mathbf{e}_k$ denoting the $k$-th canonical vector. In \eqref{eq:P1}, sparsity is imposed on $\mathbf{X}$ using \emph{elastic net} regularization \cite{elasticNet_Hastie_2005} with parameters $\mu,\lambda>0$; the elastic net regularization tends to be preferred to the plain $\ell_1$  regularization (a.k.a. LASSO), since it better preserves group patterns in the variables,   
especially when there are highly correlated variables. 
Note that Problem \eqref{eq:P1} is a nonconvex optimization problem due to the bi-convex structure of the objective functions.
The lack of convexity has motivated a lot of interest in pursuing  approximate solutions that approach optimal performance at moderate complexity; some  recent efforts are documented in   \cite{Aharon_KSVD_2006,AroraGM13,Zhang_DKSVD_2010,YaghoobiBD09,
Mairal_Bach_2010,NIPS2012_4571,RazaviyaynTsengLuo14,Scutari-Facchinei-SagratellaTSP15,Dan-Facch-Kung-ScutTSP15}. All the algorithms therein are centralized, i.e.,  they require the data matrix $\mathbf{S}$ to be centrally available.
   
In many large-scale signal processing and machine learning problems, data are not necessarily  centrally available, but collected/stored in multiple locations; {for example, consider} sensor,  cloud, or cluster-computer networks. In these scenarios, sharing local information with a central processor might be either unfeasible or not economical/efficient, owing to
the large size of the network/data, dynamicity of network topology, energy constraints, and/or privacy issues.
Common to all the aforementioned problems is the necessity of performing
a completely decentralized computation/optimization.

 Motivated by these observations,  we aim at developing
a  solution  method for \eqref{eq:P1}  in a   distributed setting wherein  the data matrix $\mathbf{S}$ is spread over a  network with (possibly)  time-varying topology, and each node has access \emph{only} to some portions of $\mathbf{S}$. 
We are not aware of any \emph{distributed, convergent} scheme for such class of problems; some attempts have been documented recently in \cite{Liang_ADMM_2014,ShakeriRB14,Chen_DDL_2015,Wai_EXTRA_AO_2015,Chainais13}. Although there are substantial differences between these methods, they can be generically abstracted as combinations of local (primal or dual ascent) descent steps followed by variable exchanges and averaging of information among neighboring nodes. However,   theoretical  convergence of these methods remains an open question.  This is mainly due to the fact that these schemes  exploit decomposition techniques suitable for (strongly) convex problems, whereas Problem \eqref{eq:P1} is nonconvex (and may lack zero-duality gap). Furthermore,  numerical  results therein are contradictory; for instance,  some of the aforementioned schemes are observed not to converge while some others fail to reach asymptotic agreement among the local copies of the dictionary variables.  Other relevant papers are \cite{Zhu-Martinez2,Bianchi-Jakubowicz,NEXT} wherein multiagent nonconvex optimization over networks is studied. While interesting, the approaches in  \cite{Zhu-Martinez2,Bianchi-Jakubowicz,NEXT}  are not able to handle nonconvex problems in the form of \eqref{eq:P1} (see Sec. \ref{Discussion_alg} for more details). 

The major contribution of this work is to propose the first \emph{provably convergent} algorithmic framework for the distributed DL problem over (possibly) time-varying network topologies.   The crux of the framework is a novel convexification-decomposition
technique that hinges on our recent SCA methods
 \cite{Scutari-Facchinei-SagratellaTSP15,Dan-Facch-Kung-ScutTSP15} while leveraging i) a gradient tracking mechanism   
to locally estimate the missing global information; and ii)   a consensus step  to distribute the computation as well as  propagate the needed information over the network. 
Asymptotic convergence of the proposed algorithm to stationary solutions of \eqref{eq:P1} is proved. Preliminary numerical results show that the proposed scheme compare favorably with other decentralized state-of-the-art algorithms. 

The rest of the paper is organized as follows. Sec. \ref{sec:distributed-DL} describes the  distributed DL problem along with the network setting. The design of the new algorithm is addressed in Sec. \ref{sec:alg-design}. 
Numerical results are presented in Sec. \ref{num_results} and some conclusions are drawn in  Sec. \ref{conclusions}.

\section{Problem Description}\label{sec:distributed-DL}
Consider Problem \eqref{eq:P1} defined over a network composed of $I$ autonomous agents (nodes). Each agent $i$ owns a subset of columns of the data matrix $\mathbf{S}\triangleq [\mathbf{S}_1,\ldots ,\mathbf{S}_I]$,  say   $\mathbf{S}_i\in \mathbb{R}^{M\times n_i}$,   and controls the corresponding part  of the sparse representation matrix $\mathbf{X}\triangleq [\mathbf{X}_1,\ldots ,\mathbf{X}_I]$, where  $\mathbf{X}_i\in \mathbb{R}^{K\times n_i}$ and $\sum_i n_i=N$.  
 Then, Problem \eqref{eq:P1} can be rewritten as\vspace{-.2cm}
\begin{equation}\tag{P2}\label{eq:P2}
\begin{array}{c l}
{\displaystyle \min_{\mathbf{D},\{\mathbf{X}_i\}_{i=1}^I}} & \displaystyle\sum_{i=1}^{I}\bigg[\underset{\triangleq f_{i}\left(\mathbf{D},\mathbf{X}_{i}\right)}{\underbrace{\dfrac{{1}}{2}\left\Vert \mathbf{S}_{i}-\mathbf{D}\mathbf{X}_{i}\right\Vert^2_F }}+\lambda\norm{\mathbf{X}_{i}}_{1,1}+ {\mu} \norm{\mathbf{X}_i}_F^2\bigg]
\\
\text{s.t.} & \mathbf{D}\in\mathcal{D}.
\end{array}
\end{equation}
 Note that each agent $i$ knows only  its ``local'' function  $f_i$ (along with $\lambda$, $\mu$, and the set $\cal D$). We  remark that,  even though in \eqref{eq:P2} we assumed that the data matrix $\mathbf{S}$ is partitioned by columns, 
 the  proposed algorithm can  be readily applied also to scenarios wherein   $\mathbf{S}$ is partitioned by rows (and thus the dictionary matrix is partitioned accordingly); see the journal version of this work.\smallskip 
 
\noindent \textbf{Network Topology.}
 Time is slotted and, at each time-slot $\nu$, the network of the $I$ agents is modeled as a digraph $\mathcal{G}^{\nu}=(\mathcal{V},\mathcal{E}^\nu)$, where $\mathcal{V}=\{1,\ldots,I\}$ is the vertex  set (i.e., the set of agents)  and $\mathcal{E}^\nu$ is the set of (possibly) time-varying directed edges: $(i,j)\in \mathcal{E}^\nu$ if there is a link from $j$ to $i$ (i.e., agent $j$ can send information to agent $i$) at time $\nu$.  The set of \emph{in-neighborhood} of agent $i$ at time $\nu$ is defined as  $\mathcal{N}_i^\nu=\{j:(i,j)\in\mathcal{E}^\nu\}\cup\{i\}$; it is the set of
 agents  that can communicate with node $i$ at time $\nu$.   To let  information   propagate over the network, we make the following standard assumption on the network connectivity.\smallskip
 
\noindent \textbf{Assumption A (on the network connectivity).} {\it
The sequence of graphs $\mathcal{G}^\nu$ is B-strongly connected, i.e., there exists a finite integer $B > 0$ such that the graph $\mathcal{G}^\nu=(\mathcal{V},\mathcal{E}_B^\nu)$, with  
$\mathcal{E}_B^\nu=\bigcup_{\nu=kB}^{(k+1)B-1}\mathcal{E}^\nu$, is strongly connected for all $k\geq0$.}\smallskip
 
 To the best of our knowledge, this is the weakest condition under which one can prove convergence of distributed algorithms on time-varying directed networks.

\section{Algorithmic design of D$^2$L}\label{sec:alg-design}
 The design of a distributed scheme for \eqref{eq:P2} faces two main challenges, namely: the nonconvexity of each $f_i$  and the lack of global information on all $f_i$. To cope with these issues,
we propose to combine alternating optimization  (Step 1 below) with
consensus mechanisms (Step 2), as described next.\smallskip 

\noindent \textbf{Step 1: Local Optimization.}
Each agent $i$     maintains    a ``local copy''  $\mathbf{D}_{(i)}$ of the common dictionary $\mathbf{D}$   and  controls its private variables $\mathbf{X}_i$. We denote by $(\mathbf{D}_{(i)}^\nu,\mathbf{X}_i^\nu)$ the value of such a pair at iteration $\nu$. The goal is that each agent updates   $(\mathbf{D}_{(i)}^\nu,\mathbf{X}_i^\nu)$ to $(\mathbf{D}_{(i)}^{\nu+1},\mathbf{X}_i^{\nu+1})$ towards a solution of Problem \eqref{eq:P2}. 
However, (directly) solving \eqref{eq:P2} is too costly, due to the nonconvexity of  $f_i(\mathbf{D}_{(i)},\mathbf{X}_i)$, and it is not even feasible (because of the lack of global knowledge of problem, i.e., agent $i$ does not have access to $f_j,\forall j\neq i$). Therefore, the idea is to somehow approximate   \eqref{eq:P2} so that   each agent  computes the new iteration  \emph{locally} and \emph{efficiently}. Since each $f_i$ is not jointly convex  in  $(\mathbf{D}_{(i)},\mathbf{X}_i)$ but bi-convex (i.e., convex in $\mathbf{D}_{(i)}$ and $\mathbf{X}_{i}$, \emph{separately}),    a  natural approach is to update $\mathbf{D}_{(i)}$ and $\mathbf{X}_i$ in an alternating fashion. Specifically, fixing  $\mathbf{X}_i=\mathbf{X}_i^\nu$,  $\mathbf{D}_{(i)}$ is updated solving the following strongly convex
approximation of \eqref{eq:P2}:\vspace{-0.1cm}
\begin{equation}
\label{D_hat_subproblem}
\hspace{-0.2cm}\mathbf{\widehat{D}}_{(i)}^\nu\overset{\Delta}{=}\underset{\mathbf{D}_{(i)}\in\mathcal{D}}
{\mathrm{argmin}}~\tilde{f}_i(\mathbf{D}_{(i)};\mathbf{D}_{(i)}^\nu,\mathbf{X}_i^\nu)+
\left\langle\boldsymbol{\Pi}_i^\nu,\mathbf{D}_{(i)}-\mathbf{D}_{(i)}^\nu\right\rangle,\vspace{-0.1cm}
\end{equation}
where $\left\langle\mathbf{A},\mathbf{B}\right\rangle\triangleq \text{tr}(\mathbf{A}^T\mathbf{B})$ and  $\tilde{f}_i(\mathbf{D}_{(i)};\mathbf{D}_{(i)}^\nu,\mathbf{X}_i^\nu)$ is defined as\vspace{-0.2cm}
\begin{equation}
\label{eq:ftilde}
\hspace{-0.2cm}\tilde{f}_i(\mathbf{D}_{(i)};\mathbf{D}_{(i)}^\nu,\mathbf{X}_i^\nu)\triangleq f_i(\mathbf{D}_{(i)},
\mathbf{X}_i^\nu)+\frac{\tau_{D,i}^\nu}{2}\, ||\mathbf{D}_{(i)}-\mathbf{D}_{(i)}^\nu||^2_F,
\end{equation}
with $\tau_{D,i}^\nu>0$ and $\boldsymbol{\Pi}_i^{\nu}\triangleq \sum_{j\neq i}\nabla_D f_j(\mathbf{D}_{(i)}^\nu,\mathbf{X}_j^\nu)$
is the \textit{linearization}  of the \emph{unknown} term
$\sum_{j\neq i}f_j(\mathbf{D}_{(i)}^\nu,\mathbf{X}_j^\nu)$.
Note that the quadratic term
$\frac{\tau_{D,i}^\nu}{2} ||\mathbf{D}_{(i)}-\mathbf{D}_{(i)}^\nu||_F^2$ in \eqref{eq:ftilde} serves to make $\tilde f_i$ strongly convex, so that \eqref{D_hat_subproblem} has a unique solution, denoted by $\mathbf{\widehat{D}}_{(i)}^\nu$. The direct use of $\widehat{\mathbf{D}}_{(i)}^\nu$ as the new local estimate $\mathbf{{D}}_{(i)}^{\nu+1}$ would not help to establish convergence  for two reasons: (i) $\mathbf{\widehat{D}}_{(i)}^\nu$ might be a too ``aggressive'' update   and (ii) we have not introduced any mechanism yet to ensure that the local copies $\mathbf{\widehat{D}}_{(i)}^\nu$ eventually agree among all agents. To cope with these  two issues,  we introduce a step-size  in the update of the dictionary:\vspace{-.1cm}
\begin{equation}
\label{u_update_ver0}
\mathbf{U}_{(i)}^\nu\triangleq\mathbf{D}_{(i)}^\nu+\gamma^\nu(\widehat{\mathbf{D}}_{(i)}^\nu- \mathbf{D}_{(i)}^\nu),\vspace{-.1cm}
\end{equation}
where $\gamma^\nu$ is a positive scalar to be  properly chosen, see Assumption C in Sec.\,\ref{Discussion_alg}).

We now consider the update of private variables  $\mathbf{X}_i$. Fixing $\mathbf{D}_{(i)}=\mathbf{U}_{(i)}^\nu$, 
agent $i$ computes the new update $\mathbf{X}_i^{\nu+1}$ solving the following strongly convex optimization problem
\begin{equation}\label{eq:Xupdate}
\mathbf{X}_i^{\nu+1}\triangleq \underset{\mathbf{X}_i\in\mathbb{R}^{K\times n_i}}{\text{argmin}} \tilde{h}_i(\mathbf{X}_i;
\mathbf{U}_{(i)}^\nu,\mathbf{X}_i^\nu)+\lambda\norm{\mathbf{X}_i}_{1,1}+ {\mu} \norm{\mathbf{X}_i}^2_F,\vspace{-0.2cm}
\end{equation}
where
\begin{equation}
\label{eq:htilde}
\tilde{h}_i(\mathbf{X}_i;
\mathbf{U}_{(i)}^\nu,\mathbf{X}_i^\nu)\triangleq f_i(\mathbf{U}_{(i)}^\nu,
\mathbf{X}_i)+ \frac{\tau_{X,i}^\nu}{2}\norm{\mathbf{X}_i-\mathbf{X}_i^\nu}_F^2,
\end{equation}
and  $\tau_{X,i}^\nu$ is a positive scalar (to be properly chosen together with $\tau_{D,i}^\nu$, see Assumption C in Sec.\,\ref{Discussion_alg}).\smallskip 

\noindent \textbf{Step 2: Consensus.}
To force the asymptotic agreement among the $\mathbf{D}^\nu_{(i)}$'s,  a consensus-based step is employed on $\mathbf{U}_{(i)}^\nu$'s.  Each agent $i$ computes the new update $\mathbf{D}_{(i)}^{\nu+1}$ as
\begin{equation}
\label{consensus_on_D}
\textstyle
\mathbf{D}_{(i)}^{\nu+1}\triangleq\displaystyle\sum_{j\in\mathcal{N}_i^\nu} w_{ij}^\nu\mathbf{U}_{(j)}^\nu
\end{equation}
where   $( w_{ij}^\nu)_{i,j=1}^I$ is a set of   weights matching the network topology $\mathcal G^\nu$ at time slot $\nu$, in the sense defined next.  

\noindent \textbf{Assumption B (on the weight matrix).} {\it The weights  $\mathbf{W}^\nu\triangleq (w_{ij}^\nu)_{i,j=1}^I$  satisfy the following conditions: i) for every  $\nu\geq 0$, 
\begin{align}
\label{weights}
w_{ij}^\nu=\left\{
             \begin{array}{ll}
               \theta\in[\vartheta,1] & \mathrm{if}~ j\in \mathcal{N}_i^\nu, \\
               0 & \mathrm{otherwise,}
             \end{array}
           \right.
\end{align}
for some $\vartheta\in (0,1)$; ii) $\mathbf{W}^\nu\,\mathbf{1}=\mathbf{1}$; and  iii) $ \mathbf{1}^T \mathbf{W}^\nu=\mathbf{1}^T$.}

\noindent \textbf{Toward a fully distributed implementation.} The computation of  $\widehat{\mathbf{D}}_{(i)}^\nu $, and consequently, the update of $\mathbf{D}^{\nu+1}_{(i)}$ in \eqref{consensus_on_D} have a severe drawback: the evaluation of $\boldsymbol{\Pi}_i^\nu$ in  \eqref{D_hat_subproblem} would require the knowledge of $\nabla_D f_j(\mathbf{D}^\nu_{(i)}, \mathbf{X}_j^\nu)$ for all $j\neq i$, which is not available  to agent $i$.
To cope with this issue, we replace $\boldsymbol{\Pi}_i^\nu$ in \eqref{D_hat_subproblem}  with a ``local estimate'', denoted by $\boldsymbol{\widetilde{\Pi}}_{i}^\nu$, and solve instead,
\begin{equation}
\label{D_tilde_subproblem}
\mathbf{\widetilde{D}}_{(i)}^\nu\overset{\Delta}{=}\underset{\mathbf{D}_{(i)}\in\mathcal{D}} {\text{argmin}}\quad\tilde{f}_i(\mathbf{D}_{(i)};\mathbf{D}_{(i)}^\nu,\mathbf{X}_i^\nu)+\left\langle
\boldsymbol{\widetilde{\Pi}}_{i}^\nu,\mathbf{D}_{(i)}-\mathbf{D}_{(i)}^\nu\right\rangle. 
\end{equation}
The question now becomes how to update   $\boldsymbol{\widetilde{\Pi}}_{i}^\nu$ using only \emph{local} information so that $\boldsymbol{\widetilde{\Pi}}_{i}^\nu$ will track the right $\boldsymbol{{\Pi}}_{i}^\nu$. 
  Rewriting $\boldsymbol{\Pi}_i^\nu$ as\vspace{-0.2cm}
\begin{equation}
\label{Pi_update_idea}
\boldsymbol{\Pi}_i^\nu=I \cdot\underbrace{\left(\textstyle\dfrac{1}{I}\displaystyle \sum_{j=1}^I\nabla_D f_j(\mathbf{D}^\nu_{(i)},
\mathbf{X}^\nu_j)\right)}_{\overset{\Delta}{=}\overline{\boldsymbol{\Theta}}_i^\nu}-\nabla_D f_i(\mathbf{D}^\nu_{(i)}, \mathbf{X}^\nu_i),
\end{equation}
we propose   to compute $\boldsymbol{\widetilde{\Pi}}_i^{\nu}$ mimicking \eqref{Pi_update_idea}:
\begin{equation}
\label{Pi_update}
\boldsymbol{\widetilde{\Pi}}_i^\nu=I\cdot\boldsymbol{{\Theta}}_i^\nu-\nabla_D f_i(\mathbf{D}^\nu_{(i)},\mathbf{X}^\nu_i),
\end{equation}
\vskip-0.1cm
\noindent
where $\boldsymbol{{\Theta}}_i^\nu$ is a local auxiliary variable (controlled by user $i$) whose task is to asymptotically track $\overline{\boldsymbol{\Theta}}_i^\nu$. This can be done by leveraging the  tracking mechanism first  introduced in \cite{NEXT}:
\begin{align}
\label{theta_update}
\boldsymbol{{\Theta}}_i^{\nu+1}= &\sum_{j\in\mathcal{N}_i^\nu} w_{ij}^\nu \boldsymbol{{\Theta}}_j^{\nu}
\\
&+\nabla_D f_i(\mathbf{D}_{(i)}^{\nu+1},\mathbf{X}_i^{\nu+1})-\nabla_D f_i(\mathbf{D}_{(i)}^{\nu},
\mathbf{X}_i^{\nu})\nonumber
\end{align}
with $\boldsymbol{{\Theta}}_i^0\overset{\Delta}{=}\nabla_D f_i(\mathbf{D}_{(i)}^0,\mathbf{X}_i^0)$ for all $i=1,2,\ldots,I$.

Because of  above modifications towards a distributed implementation, the update rule \eqref{u_update_ver0} is impacted and needs to be properly modified by replacing $\widehat{\mathbf{D}}_{(i)}^\nu$ with $\widetilde{\mathbf{D}}_{(i)}^\nu$, which reads\vspace{-0.15cm}
\begin{equation}
\label{u_update}
\mathbf{U}_{(i)}^\nu=\mathbf{D}_{(i)}^\nu+\gamma^\nu(\mathbf{\widetilde{D}}_{(i)}^\nu- \mathbf{D}_{(i)}^\nu).
\end{equation}
We can now formally introduce the proposed Distributed Dictionary Learning (D$^2$L) algorithm, as given in Algorithm \ref{alg1}. Its convergence
properties   are stated in Theorem \ref{th:conver}.

\begin{algorithm}[t]
$\textbf{Initialize}:$ $\mathbf{X}_i^0=\mathbf{0},~\mathbf{D}_{(i)}^0\in\mathcal{D},~\boldsymbol{{\Theta}}^0_i=\nabla_D f_i(\mathbf{D}^0_{(i)},\mathbf{X}^0_i),$ 

\quad\quad\qquad~ $\boldsymbol{\widetilde{\Pi}}_i^0=I\cdot\boldsymbol{{\Theta}}_i^0-\nabla_D f_i(\mathbf{D}_{(i)}^0,\mathbf{X}_i^0),~\forall i$; set  $\nu=0;$
\\
\texttt{S1.} If $(\mathbf{D}_{(i)}^\nu,\mathbf{X}_i^\nu)$ satisfies stopping criterion for all $i$'s: \texttt{STOP};\smallskip
\\
\texttt{S2.} \textbf{Local Updates:} Each agent $i$ computes:
\begin{enumerate}[label=(\alph*)]
\item $\widetilde{\mathbf{D}}_{(i)}^\nu$ and $\mathbf{U}_{(i)}^\nu$ according to \eqref{D_tilde_subproblem} and \eqref{u_update};
\item $\mathbf{X}_i^{\nu+1}$ according to \eqref{eq:Xupdate};
\end{enumerate}
\texttt{S3.} \textbf{Broadcasting:} Each agent $i$ collects data
from its current neighbors and updates:
\begin{enumerate}[label=(\alph*)]
\item $\mathbf{D}_{(i)}^{\nu+1}$ according to \eqref{consensus_on_D};
\item $\boldsymbol{{\Theta}}_i^{\nu+1}$ and $\boldsymbol{\widetilde{\Pi}}_i^{\nu+1}$ according to  \eqref{theta_update} and \eqref{Pi_update};
\end{enumerate}
\texttt{S4.} Set $\nu+1\to\nu$, and go to \texttt{S1}. 
\protect\caption{\textbf{: Distributed Dictionary Learning (D$^2$L)}}
\label{alg1} 
\end{algorithm}

\subsection{Convergence of Algorithm \ref{th:conver}} 
\label{Discussion_alg}
In Algorithm\,1, there are some parameters to be tuned, namely:   i) the step-size $\gamma^\nu$; and ii)  the proximal coefficients $(\tau_{X,i}^\nu)_{i=1}^I$ and $(\tau_{D,i}^\nu)_{i=1}^I$. While several choices are possible for these quantities, some minimal  conditions need to be satisfied to guarantee convergence of Algorithm\,1 as well as asymptotic consensus. More specifically, we need the following.\smallskip 

\noindent \textbf{Assumption C (on the free parameters).} 
The parameters  $\gamma^\nu$, $(\tau_{X,i}^\nu)_{i=1}^I$ and $(\tau_{D,i}^\nu)_{i=1}^I$ are chosen such that:
\begin{description}
\item[\textbf{(C1)}]  $\{\gamma^\nu\}_\nu$  satisfies:  $\gamma^\nu\in[0,1]$ for all $\nu\geq 1$,  
$\sum_{\nu=1}^\infty \gamma^\nu=\infty$, and $\sum_{\nu=1}^\infty\left(\gamma^\nu\right)^2<\infty$;
\item[\textbf{(C2)}] The sequences $\{\tau_{X,i}^\nu\}_{\nu}$ and $\{\tau_{D,i}^\nu\}_{\nu}$ are bounded and uniformly positive, for all $i=1,\ldots, I$. Additionally, $\tau_{X,i}^\nu=\max(\epsilon,$ $\sigma_{\mathrm{max}}(\mathbf{U}_{(i)}^\nu)^2)$, where $\epsilon>0$ is  arbitrary, and $\sigma_{\mathrm{max}}(\mathbf{U}_{(i)}^\nu)$ is the maximum singular value of $\mathbf{U}_{(i)}^\nu$.\smallskip
\end{description}
 
We can now provide the main convergence results, as stated in   Theorem\,1 below.  The proof of the  theorem is quite involved and   omitted here for lack of space; see the journal version of this work.
\smallskip 
\begin{theorem}
\label{th:conver}
\it
Let $\{(\mathbf{D}^\nu_{(i)},\mathbf{X}^\nu_i)_{i=1}^I\}_\nu$ be the sequence generated by Algorithm\,1  and let $\overline{\mathbf{D}}^\nu\triangleq \frac{1}{I}\sum_{i=1}^I\mathbf{D}_{(i)}^\nu$. Suppose that Assumptions A-C are satisfied, then the following holds: \begin{itemize}\item[(i)]  $\{(\overline{\mathbf{D}}^\nu,\mathbf{X}^\nu)\}_\nu$ is bounded and any of its limit points is a stationary solution of Problem \eqref{eq:P2} \emph{[}and thus \eqref{eq:P1}\emph{]}; 
and \item[(ii)] all  $\{\mathbf{D}_{(i)}^\nu\}_\nu$ asymptotically reach consensus, i.e.,  $
\lim_{\nu\rightarrow\infty}||\mathbf{D}_{(i)}^\nu-\overline{\mathbf{D}}^\nu||=0$ for all  $i=1,2,\ldots,I$.\end{itemize}
\end{theorem}

Roughly speaking, Theorem 1 states two main results: 1) (subsequence) convergence of $\{(\overline{\mathbf{D}}^\nu,\mathbf{X}^\nu)\}_{\nu}$  to a stationary solution of \eqref{eq:P1}; and 2)  asymptotic  agreement of all  $\mathbf{D}_{(i)}^\nu$  on the limit point of $\{\overline{\mathbf{D}}^\nu\}_{\nu}$.\vspace{-0.1cm}

\subsection{Discussion}
\label{Discussion_alg}\vspace{-0.1cm}
\noindent \textbf{Convergence:} To the best of our knowledge, Algorithm\,1 is the  first   distributed algorithm   for the DL problem \eqref{eq:P1}, with  convergence guarantees. Our results can be contrasted with  \cite{Liang_ADMM_2014,ShakeriRB14,Chen_DDL_2015,Wai_EXTRA_AO_2015,Chainais13} wherein gradient schemes tailored with consensus/diffusion updates  are employed for some instance of \eqref{eq:P1}.  The aforementioned schemes \emph{do not have any convergence guarantees}: it is postulated that the sequence generated by the algorithms is convergent   (see, e.g., \cite{Wai_EXTRA_AO_2015,Chainais13}),   and then concluded {that any limit point is a stationary
solution of the problem. Furthermore,   some of these schemes do not even achieve  consensus among the local variables.   Finally, we remark that our previous results  \cite{NEXT} are not applicable to \eqref{eq:P1}, since:   i)  \cite{NEXT} cannot handle private variables, i.e., $\mathbf{X}_i$'s; and ii) convergence of \cite{NEXT} (and   of most of distributed gradient schemes in the literature \cite{Liang_ADMM_2014,ShakeriRB14,Chen_DDL_2015,Wai_EXTRA_AO_2015,Chainais13})  require some technical properties that are not satisfied by the functions in  \eqref{eq:P2} [e.g, boundedness and  Lipschitzianity of the gradient of $f_i$].

\noindent \textbf{Solutions of the subproblems:} The update of the variables $(\mathbf{{D}}^\nu_{(i)},\mathbf{X}^\nu_i)$ calls for the solutions of strongly convex problems [cf. \eqref{eq:Xupdate} and \eqref{D_tilde_subproblem}]. One can of course rely on standard solvers for convex problems. 
To alleviate the computational burden, one can alternatively choose   different surrogates  $\tilde{f}_i$s in \eqref{eq:Xupdate} and \eqref{D_tilde_subproblem} that deliver closed form solutions for $\widetilde{\mathbf{D}}_{(i)}^\nu$ and ${\mathbf{X}}_{(i)}^{\nu+1}$. More specifically, considering \eqref{D_tilde_subproblem}, one can linearize $f_i$, that is,
\begin{align}
\label{eq:ftilde2}
\tilde{f}_i(\mathbf{D}_{(i)};\mathbf{D}_{(i)}^\nu,\mathbf{X}_i^\nu)=&\left\langle \nabla_D f_i(\mathbf{D}^\nu_{(i)}, \mathbf{X}_i^\nu),\mathbf{D}_{(i)}-\mathbf{D}_{(i)}^\nu\right\rangle
\\
&\qquad\qquad\quad~+\frac{\tau_{D,i}^\nu}{2}\, \big\|\mathbf{D}_{(i)}-\mathbf{D}_{(i)}^\nu\big\|^2,\nonumber
\end{align}
which  leads to the following closed form solution for $\widetilde{\mathbf{D}}_{(i)}^\nu$:
\begin{equation} 
\label{closed_form_dictionary_lin}
\widetilde{\mathbf{D}}_{(i)}^\nu=P_\mathcal{D}\left[\mathbf{D}_{(i)}^\nu-\frac{1}{\tau_{D,i}^\nu}\left(\nabla_{D} f_i(\mathbf{D}_{(i)}^\nu,\mathbf{X}_i^\nu)+\boldsymbol{\widetilde{\Pi}}_i^\nu\right)\right].
\end{equation}

Consider now the sparse coding subproblem \eqref{eq:Xupdate}. If  $\tilde{h}_i$ is chosen as in \eqref{eq:htilde}, the update of the local variables  $\mathbf{X}_i^{\nu+1}$
 reduces to solving a LASSO problem (see, e.g.,  \cite{Scutari-Facchinei-SagratellaTSP15,Dan-Facch-Kung-ScutTSP15} for recent
efficient algorithms for large-scale LASSO problems). To avoid to solve
a LASSO problem, we can alternatively use the linearization of $f_i$ as $\tilde{h}_i$, that is,\vspace{-0.2cm}
\begin{align}
\label{eq:htilde2}
\tilde h_i(\mathbf{X}_i;\mathbf{U}_{(i)}^\nu,\mathbf{X}_i^\nu)&= \left\langle\nabla_{X_i} f_i(\mathbf{U}_{(i)}^\nu,\mathbf{X}_i^\nu),\mathbf{X}_i-\mathbf{X}_i^\nu\right\rangle
\\
&\qquad\qquad\qquad\qquad+\frac{\tau_{X,i}^\nu}{2}\norm{\mathbf{X}_i-\mathbf{X}_i^\nu}^2,\nonumber
\end{align}
which leads to the following closed form solution for $\mathbf{X}_i^{\nu+1}$: \begin{equation}
\mathbf{X}_i^{\nu+1}
=\frac{\tau^\nu_{X,i}}{2\,\mu+\tau^\nu_{X,i}}\cdot\mathcal{T}_{\frac{\lambda}{\tau^\nu_{X,i}}} \left(\mathbf{X}^\nu_i-\frac{1}{\tau^\nu_{X,i}}\nabla_{X_i} f_i(\mathbf{U}_{(i)}^\nu,\mathbf{X}_i^\nu) \right),
\label{closed_form_X_linearization}
\end{equation}
where  $\mathcal T_\theta(x)\triangleq \max(|x|-\theta,0)\cdot \text{sign}(x)$ is the soft-thresholding operator and is applied element-wise in \eqref{closed_form_X_linearization}. In Sec.\,\ref{num_results}, we present some  numerical results comparing the two versions of Algorithm\,1, based on the  choices of $\tilde{h}_i$ in \eqref{eq:htilde} and \eqref{eq:htilde2}, respectively.  We remark that the convergence results stated in Theorem \ref{th:conver} remain valid  also for the aforementioned  new choices of surrogate functions;  see the journal version of this work.

\noindent\textbf{Tuning of free parameters:} There are three set of  parameters to tune in Algorithm\,1, namely: i) the step-size $\gamma^\nu$; ii) the proximal coefficients $(\tau_{X,i}^\nu)_{i=1}^I$ and $(\tau_{D,i}^\nu)_{i=1}^I$; and iii) the weights $(w_{ij}^\nu)_{i,j=1}^I$. Theorem\,\ref{th:conver} offers some flexibility in the choice of these parameters (cf.\,Assumptions B and C). For instance, the condition on the step-size--Assumption C1--ensures that the sequence decays to zero, but not too fast. There are many diminishing step-size rules in the literature satisfying  this condition; see, e.g., \cite{Bertsekas2}. An effective  step-size rule that we used in our experiments (see \cite{Scutari-Facchinei-SagratellaTSP15,Dan-Facch-Kung-ScutTSP15}) is $\gamma^\nu=\gamma^{\nu-1}(1-\epsilon \gamma^{\nu-1})$ with $\gamma^0\in (1,2]$ and $\epsilon\in (0,1/\gamma^0)$.  
The proximal coefficients can be set as $\tau^\nu_{D,i}=\tilde{\epsilon}$, for all $i$ and $\nu$, where $\tilde{\epsilon}>0$; and the explicit expression for $\tau^\nu_{X,i}$ is given in Assumption\,C2. 
Note that the above choices of step-size and proximal coefficients do not require any form of centralized coordination among the agents, which is a key feature in our
distributed environment.
Finally, referring to the weights $(w_{ij}^\nu)_{i,j=1}^I$, several choices satisfying Assumption B are available, see  \cite{NEXT} and references therein for some examples.  


\section{Numerical Results}\label{num_results} 
In this section we present some numerical results comparing
  Algorithm\,1 with the (distributed) ATC algorithm \cite{Chainais13}. For Algorithm\,1  we simulated  two instances, namely: i) one based on the surrogates    \eqref{eq:ftilde2} and \eqref{eq:htilde}, which we will refer to   as  ``Plain D$^2$L'';  and ii) one using the surrogates \eqref{eq:ftilde2} and  \eqref{eq:htilde2}, which will be termed  ``Linearized D$^2$L''.  

\noindent \textbf{Setting and tuning:} We consider denosing a $512\times 512$ pixels corrupted boat image in a distributed setting.  The data set $\mathbf{S}$ is composed of the stacked $8\times 8$ sliding patches of the image. The size of the dictionary and the sparse representation matrices $\mathbf{X}_i$ are $64\times 64$ and $64\times 255,150$, respectively (overall, the number of variables is around $16$ million), and the parameters in \eqref{eq:P2} are set to $2\,\mu=\lambda=1/8$ and $\alpha=1$. We simulated a time-invariant undirected connected network composed of 150 agents. For all the algorithms, the local copies $\mathbf{D}_{(i)}$'s are initialized to random patches of the local data and  $\mathbf{X}_i$'s are initialized to zero. In both versions of our algorithms, the diminishing step size sequence $\gamma_\nu$ is generated according to $\gamma^\nu=\gamma^{\nu-1}(1-\epsilon \gamma^{\nu-1})$ with  $\gamma^0=0.5$ and $\epsilon=0.1$. 
The weights $(w_{ij}^\nu)_{i,j=1}^I$ in the consensus steps are chosen according to the Metropolis rule \cite{Xiao}.

\noindent\textbf{Merit Function:}
We introduce the following functions to measure    progresses of the  algorithms towards stationarity of \eqref{eq:P1} and attainment of consensus. Using \cite{Scutari-Facchinei-SagratellaTSP15}, it is not difficult to check that $\Delta^\nu=||\mathrm{vec}(\boldsymbol{\Delta}_D^\nu,\boldsymbol{\Delta}_X^\nu)||_\infty$ is a valid distance from stationarity, where 
\begin{equation}
\begin{aligned}
\boldsymbol{\Delta}_D^\nu=\overline{\mathbf{D}}^\nu-\widehat{\mathbf{D}}^\nu,
\quad\boldsymbol{\Delta}_X^\nu=\mathbf{X}^\nu-\widehat{\mathbf{X}}^\nu
\end{aligned}
\end{equation}
and \vspace{-0.3cm}
\begin{equation}
\label{stat_merit_subproblems}
\begin{aligned}
&\qquad\qquad\qquad\widehat{\mathbf{D}}^\nu\triangleq \underset{\mathbf{D}\in\mathcal{D}}
{\mathrm{argmin}}~\sum_{i=1}^I\tilde{f}_i(\mathbf{D};\overline{\mathbf{D}}^\nu,\mathbf{X}_i^\nu),
\\
&\widehat{\mathbf{X}}^\nu\triangleq \underset{\{\mathbf{X}_i\}_{i=1}^I}
{\mathrm{argmin}}~\sum_{i=1}^I\tilde{h}_i(\mathbf{X}_i;\overline{\mathbf{D}}^\nu,\mathbf{X}_i^\nu)+\lambda\norm{\mathbf{X}_i}_{1,1}+\mu\norm{\mathbf{X}_i}_F^2,
\end{aligned}
\end{equation}
where $\tilde{f}_i$ and $\tilde{h}_i$ are defined in \eqref{eq:ftilde2} and \eqref{eq:htilde2}, respectively; and   $\tau_{D,i}^\nu=1$ and $\tau_{X,i}^\nu=1$, for all $i$ and $\nu$. Note that  $\Delta^\nu$ is a continuous function of $(\overline{\mathbf{D}}^\nu,\mathbf{X}^\nu)$  and is zero if and only if the argument is a  stationary solution of \eqref{eq:P1}. The consensus disagreement is instead   
measured by computing the consensus error 
$e^\nu=\max_i||\mathrm{vec}(\mathbf{D}_{(i)}^\nu-\overline{\mathbf{D}}^\nu)||_\infty$.

In Fig. \ref{fig:obj_consErr_stationarity} we plot, for each algorithm, the above merit functions [and the objective function of (P1) evaluated at $(\overline{\mathbf{D}}^\nu, \mathbf{X}^\nu)$]   versus the number of agents' message exchanges. For ATC,  the number of message exchanges coincides with the iterations $\nu$, while for our schemes, it is $2\cdot \nu$. The figures clearly show that  both versions of Algorithm\,1 are much faster than ATC while being also guaranteed to converge (or, equivalently, they require less information exchanges than ATC). Note that  the computational cost per iteration of Plain D$^2$L  is comparable with that of ATC  (both require to solve a LASSO), whereas that of Linearized D$^2$L is cheaper.  Note also that ATC does not seem  to  reach a consensus on the local copies of the dictionary, whereas for our D$^2$L schemes consensus is reached quite early and then maintained. 

To measure the quality  of the reconstruction, we report in Table \ref{image_quality_table} the PSNR and the MSE values of the  reconstructed images achieved by the three algorithms, after $200$ and $1000$ message exchanges.\vspace{-0.1cm}
\begin{table}[H]
 \renewcommand{\arraystretch}{.8}
 \begin{tabular}{c l l l} 
\rowcolor{gray!15}
\cellcolor{gray!0}  & \parbox{1.6cm}{\center Linearized D$^2$L\\~} & \qquad\parbox{1.2cm}{\center Plain D$^2$L\\~} & \parbox{2cm}{\center ATC\\~} \\ [0.5ex] 
\cellcolor{gray!15}  \parbox{1.5cm}{\center 200 message exchanges\\~} & \parbox{1.2cm}{PSNR=27.28db\\MSE=121.4} & \parbox{1cm}{PSNR=27.32db\\MSE=120.2} & \parbox{1cm}{PSNR=26.48db\\MSE=146.2}
\\
\cellcolor{gray!15}  \parbox{1.7cm}{\center 1000 message exchanges\\~} & \parbox{1.2cm}{PSNR=27.53db\\MSE=114.6} & \parbox{1cm}{PSNR=27.65db\\MSE=111.69} & \parbox{1cm}{PSNR=27.29db\\MSE=121.23} \\  [1ex] 
\end{tabular}
\caption{Reconstructed image quality  from a noisy image, with PSNR=$20.34$db and MSE=$601.1$. }\vspace{-0.2cm}
\label{image_quality_table}
\end{table} 
The comparisons above shows that both D$^2$L schemes  attain good quality solutions already after 200 message exchanges, while ATC lags significantly behind. The superior
  performance of the  D$^2$L schemes seems  mainly due to the employed  \emph{gradient tracking} mechanism---instead of neglecting
  $\sum_{j\neq i}f_j$ in (P2) (like in ATC), each agent $i$ tracks the 
  gradient $\sum_{j\neq i}\nabla f_j$. 
  \begin{figure}[t]
\hspace{-.5cm}\includegraphics[scale=0.2]{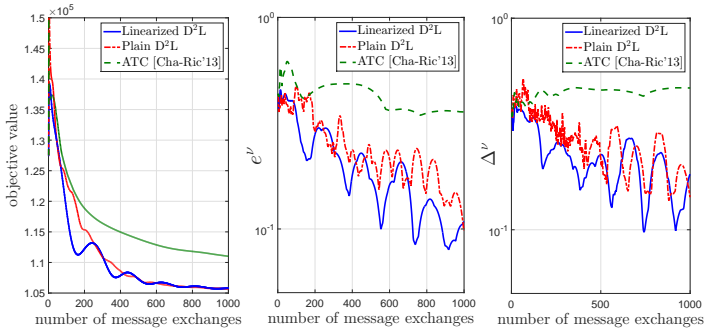}
\caption{ Objective function  (left), consensus disagreement (center), and  distance from stationarity (right) vs. number of message exchanges.}
\label{fig:obj_consErr_stationarity}\vspace{-0.6cm}
\end{figure}

\section{conclusions}\label{conclusions} 
The paper studied the distributed DL problem over (possibly) time-varying networks. We proposed the first distributed algorithmic framework  with provable convergence for this class of problems. Preliminary numerical results show promising performance for the proposed scheme. 

\bibliographystyle{ieeetr}
\bibliography{NewRefs,scutari_refs,Surbib}
\end{document}